\documentclass[10pt]{article}
\usepackage{amssymb,latexsym, amsmath}
\usepackage{amscd}

\newtheorem{theorem}{Theorem}

\newtheorem{corollary}{Corollary}

\newtheorem{lemma}{Lemma}

\newcommand{\R}{\mathbb{R}}

\newcommand{\ddim}{\mathrm{ddim\;}}
\newcommand{\dind}{\mathrm{dind\;}}
\newcommand{\corank}{\mathrm{corank\;}}
\newcommand{\g}{\mathfrak{g}}

\newcommand{\ad}{\mathrm{ad}}
\DeclareMathOperator{\Ad}{\mathrm{Ad}}
\DeclareMathOperator{\Span}{\mathrm{span\,}}
\DeclareMathOperator{\diag}{\mathrm{diag}}

\DeclareMathOperator{\rank}{\mathrm{rank}}
 \DeclareMathOperator{\pr}{\mathrm{pr}}

\begin{document}

\title{Integrability of Invariant Geodesic Flows on $n$-Symmetric Spaces}

\author{ Bo\v zidar Jovanovi\'c \\ \small Mathematical Institute, Serbian Academy of Sciences and
Arts
\\\small Kneza Mihaila 36, 11000 Belgrade, Serbia}
\date{}

\maketitle

\begin{abstract} In this paper, by modifying the argument shift
method,
we prove Liouville integrability of geodesic flows of normal
metrics (invariant Einstein metrics) on the Ledger-Obata
$n$-symmetric spaces $K^n/\diag(K)$, where $K$ is a semisimple
(respectively, simple) compact Lie group.
\\

\\

\end{abstract}

\section{Introduction}

\paragraph{Invariant Geodesic Flows.}
We study integrability of $G$-invariant geodesic flows on a class
of homogeneous spaces
\begin{equation}
 \label{L-O} Q=G/H, \quad
G=\underbrace{K\times \cdots \times K}_n, \quad
H=\diag(K)=\{(g,\dots,g)\,\vert\, g\in K\},
\end{equation}
where $K$ is a compact connected semisimple Lie group. The
homogeneous space $Q$ is diffeomorphic to the direct product
$K^{n-1}$, however as a $G$-homogeneous space it is a basic
example of a $n$-symmetric Riemannian space, see Ledger and Obata
\cite{LO}.

Let $\mathfrak g=\mathfrak k^n=\mathfrak k_1\oplus\mathfrak
k_2\oplus \dots\oplus\mathfrak k_n$, $\mathfrak
h=\{(x,\dots,x)\,\vert \, x\in \mathfrak k\}$, $\mathfrak k$ be
the Lie algebras of $G$, $H$ and $K$, respectively ($\mathfrak
k_i\cong\mathfrak k$ is the $i$-th factor). For simplicity, both
negative Killing forms on $\mathfrak g$ and $\mathfrak k$ will be
denoted by $\langle\cdot,\cdot\rangle$. Let
\begin{equation}
\mathfrak g=\mathfrak h\oplus \mathfrak v, \qquad \mathfrak
v=\{x=(x_1,\dots,x_n)\in\mathfrak k^n\, \vert \, x_1+\dots+x_n=0\}
\label{v}
\end{equation}
be the orthogonal decomposition with respect to the Killing form.

The linear subspace $\mathfrak v$ can be naturally identified with
$T_{\rho(e)}Q$, where $\rho: G\to Q=G/H$ is the canonical
projection. Then $G$-invariant metrics on $Q$, via restrictions to
$T_{\rho(e)}Q \cong \mathfrak v$, are in one-to-one correspondence
with $\Ad_H$-invariant scalar  products (e.,g., see \cite{Be})
\begin{equation}\label{I-metric}
(\,\cdot\,,\,\cdot\,)_{\mathfrak v}=\langle I(\,\cdot\,),\,\cdot\,
\rangle, \quad I: \mathfrak v\to \mathfrak v, \quad
I\circ\Ad_h=I\circ\Ad_h, \quad h\in H.
\end{equation}

The negative Killing form itself defines {\it normal (or standard)
metric} $ds^2_0$ \cite{Be}. Note that, in the case when $K$ is a
simple group, the normal $G$-invariant metric on $Q$ is Einstein
(see Wang and Ziller \cite{WZ}). Besides, Nikonorov proved that,
up to the isometry and homothety, the homogeneous space $Q$ for
$n=3$ ($n\ge 4$) admits exactly (respectively, at least) two
$G$-invariant Einstein metrics \cite{Nik}.

\paragraph{General Setting.}
Let $\mathcal F$ be a collection of functions closed under the
Poisson bracket on a Poisson manifold $(M,\{\cdot,\cdot\})$ and
let $\Lambda$ be the Poisson bivector related to
$\{\cdot,\cdot\}$. Consider the linear space $\mathrm F_x \subset
T_x^* M$ spanned by differentials of functions in $\mathcal F$.
Suppose that the numbers $\dim\mathrm F_x$ and $ \dim \ker \Lambda
\vert_{\mathrm F_x}$ are constant  almost everywhere on $M$ and
denote them by $\ddim\mathcal F$ and $\dind\mathcal F$,
respectively ({\it differential dimension} and {\it differential
index} of $\mathcal F$). The set $\mathcal F$ is called {\it
complete} if: $ \ddim{\mathcal F}+\dind{\mathcal F}=\dim
M+\corank\{\cdot,\cdot\}. $ It is {\it complete at} $x\in M$ if $
\dim\mathrm F_x + \dim \ker \Lambda \vert_{\mathrm F_x}= \dim M +
\ker\Lambda$, i.e., $\mathrm F_x$ is isotropic: $ \mathrm
F^\Lambda_x \subset \mathrm F_x$, $\mathrm F_x^\Lambda=\{\xi\in
T^*_xM \,\vert\, \Lambda_x(\xi,\mathrm F_x)=0\}.$

The Hamiltonian system $\dot f=\{f,h\}$ is {\it completely
integrable in the noncommutative sense} if it possesses a complete
set of first integrals $\mathcal F$. Then (under compactness
condition) $M$ is almost everywhere foliated by $(\dind{\mathcal
F}-\corank\{\cdot,\cdot\}$)-dimensional invariant tori. As in the
Liouville theorem, the Hamiltonian flow restricted to regular
invariant tori is quasi-periodic (see \cite{N, MF2, BJ2, Zu}).
Mishchenko and Fomenko stated the conjecture that non-commutative
integrable systems are integrable in the usual commutative sense
by means of integrals that belong to the same functional class as
the original non-commutative integrals \cite{MF2}. In the analytic
case, when $\mathcal F$ is a {finite-dimensional Lie algebra}, the
conjecture has been proved by Sadetov \cite{Sa}. The conjecture is
also proved in $C^\infty$-smooth case for infinite-dimensional
algebras (see \cite{BJ2}).

Now, let $Q=G/H$ be a homogeneous space of a compact Lie group
$G$, $\Phi: T^*Q\to \g^* $ be the momentum mapping of the natural
$G$-action on $T^*Q$, $\mathcal F_1=\Phi^*(\R[\mathfrak g])$ be
the set of Noether's functions and $\mathcal F_2$ be the set of
$G$-invariant functions, polynomial in momenta. Both $\mathcal
F_1$ and $\mathcal F_2$ are Lie subalgebras of
$(C^\infty(T^*Q),\{\cdot,\cdot\})$, where $\{\cdot,\cdot\}$ is the
canonical Poisson bracket. From Noether's theorem we have
$\{\mathcal F_1,\mathcal F_2\}=0$. Also $\mathcal F_1+\mathcal
F_2$ is a {\it complete set} of functions on $T^*Q$ (see Bolsinov
and Jovanovi\' c \cite{BJ1, BJ3}).

The Hamiltonian function $H_0$ of the normal metric $ds^2_0$ is a
Casimir  function within $\mathcal F_2$, so it Poisson commute
both with $\mathcal{F}_1$ and $\mathcal{F}_2$. Thus the geodesic
flow of the normal metric  is completely integrable in the
non-commutative sense by means of analytic functions, polynomial
in momenta.

\paragraph{Integrable Pairs.}
Within the class of Noether's integrals $\mathcal F_1$, for
example by using the argument translation method \cite{MF1}, one
can always construct a complete commutative subset of function
$\mathcal F^0\subset \mathcal F_1$ ($ \ddim {\mathcal F}^0 =
\frac12\left(\ddim \mathcal F_1 + \dind \mathcal F_1\right)$).
Thus, for the case of the geodesic flow of the normal metric, the
Mishchenko-Fomenko conjecture reduces to the construction of a
complete commutative subset $\mathcal F \subset \mathcal F_2$:
\begin{equation}
\ddim {\mathcal F} = \frac12\left(\ddim \mathcal F_2 + \dind
\mathcal F_2\right). \label{C}
\end{equation}
Indeed, from the completeness of $\mathcal F_1+\mathcal F_2$ it
follows that $\mathcal F^0+\mathcal F$ is a complete commutative
set on $T^*Q$ (see \cite{BJ1, BJ3}).

If the required subset $\mathcal F\subset \mathcal F_2$ exist, we
say that $(G,H)$ is an {\it integrable pair}. In \cite{BJ3} the
conjecture is stated that all pairs $(G,H)$ are integrable. If
$(G,H)$ is a spherical pair, in particular if $G/H$ is a symmetric
space, the algebra $\mathcal F_2$ is already commutative. In this
case we need only Noether's integrals $\mathcal F_1$ to integrate
the geodesic flow (see Mishchenko \cite{Mis}, Brailov \cite{Br}
and Mikityuk \cite{Mik}).

There are several known classes of integrable pairs (see
\cite{BJ1, BJ3, MP, DGJ}) but the general problem rest still
unsolved. For a related problem on the integrability of geodesic
flows on homogeneous spaces of noncompact Lie groups see, e.g.,
\cite{Bu, MS}.

\paragraph{Results and Outline of the Paper.}
Let $Q=G/H$ be the  Ledger-Obata $n$-symmetric space \eqref{L-O}.
By using the flag of subalgebras
\begin{equation}\label{flag}
\g_1=\mathfrak k_1\subset \g_2=\mathfrak k_1 \oplus \mathfrak k_2
\subset\dots\subset \g_n=\g=\mathfrak
k_1\oplus\dots\oplus\mathfrak k_n,
\end{equation}
we modify the argument shift method to construct a complete set of
polynomials on $\g$ with respect to the usual Lie-Poisson bracket
(Theorem \ref{prva}, Section 2). It allows us to  find a complete
commutative subset of polynomials within $\mathcal F_2$ (Theorem
\ref{n=3}, Corollary \ref{potpun}, Section 2) implying:

\begin{theorem}
The geodesic flow of the normal metric on the Ledger-Obata
$n$-symmetric space \eqref{L-O} is Liouville integrable by means
of analytic integrals, polynomial in momenta.
\end{theorem}

As a corollary, the complete commutative integrability of the
geodesic flows of invariant Einstein metrics constructed by
Nikonorov \cite{Nik} (Corollary 3, Section 3) is obtained.

\section{Liouville Integrability of Geodesic Flows}

\paragraph{$H$-invariant Euler Equations.}
Consider the left trivialization $T^*G \cong \mathfrak g \times
G$, where the identification $\g^*\cong\mathfrak g$ is given by
$\langle \cdot,\cdot\rangle$. Let $\hat I: \mathfrak g\to
\mathfrak g$ be a positive definite operator which defines
left-invariant metric $ds^2_{\hat I}$ on $G$.

The left $G$-reduction of the geodesic flow of the metric
$ds^2_{\hat I}$ is described by the Euler equations on $\mathfrak
g^* \cong\mathfrak g$:
\begin{equation}
\dot x_i=[x_i,\xi_i], \quad \xi_i=\nabla_{x_i} \hat
h(x_1,\dots,x_n)=\pr_{\mathfrak k_i} \hat A(x_1,\dots,x_n), \quad
i=1,\dots,n, \label{euler}\end{equation} where $\hat
h=\frac12\langle \hat A(x),x\rangle$ is the Hamiltonian, $\hat
A=\hat I^{-1}$ and $\pr_{\mathfrak k_i}$ is the projection to
$i$-th factor: $\pr_{\mathfrak k_i}(x_1,\dots,x_n)=x_i.$

The Euler equations are Hamiltonian with respect to the
Lie-Poisson bracket (the product of the Lie-Poisson brackets on
factors $\mathfrak k_i$):
\begin{equation}
\{f,g\}(x_1,\dots,x_n)=-\sum_{i=1}^n \langle x_i, [\nabla_{x_i}
f,\nabla_{x_i} g]\rangle. \label{LP}
\end{equation}

The  right $H$-action on $T^*G$ is Hamiltonian with momentum
mapping, in the left-trivialization, given by
\begin{equation}
\mu(x)=x_1+\dots+x_n. \label{mu}\end{equation}

The geodesic flow is invariant with respect to the right
$H$-action if and only if the Hamiltonian $\hat h$ is
$\Ad_H$-invariant, i.e,
\begin{equation*}
\langle [x,\mathfrak h], \hat A(x) \rangle=0
\,\,\Leftrightarrow\,\, \pr_{\mathfrak h} [\hat A(x),x]=0
\,\,\Leftrightarrow\,\, \sum_{i=1}^n {[\pr_{\mathfrak k_i} \hat A
(x),x_i]=0}, \label{i1}
\end{equation*}
where we used
\begin{equation}\label{proj-h}
 \pr_{\mathfrak
h}(x_1,\dots,x_n)=\frac1n(x_1+\dots+x_n,\dots,x_1+\dots+x_n).
\end{equation}

If the Hamiltonian $\hat h$ is $\Ad_H$-invariant, then the
momentum $\mu$ is preserved by geodesic flow and we can perform
the symplectic reduction of the flow to $\mu^{-1}(0)/H\cong T^*Q$.
The reduced flow is the geodesic flow of a $G$-invariant
submersion metric on $Q$.

Contrary, for a given $\Ad_H$-invariant positive definite operator
$I:\mathfrak v\to \mathfrak v$, let $ds^2_I$ be a $G$-invariant
metric defined by \eqref{I-metric}.  It can be seen as a
submersion metric of an appropriate left $G$-invariant and right
$H$-invariant metric $ds^2_{\hat I}$, simply by taking
$$
\hat I(x)=s\cdot \pr_\mathfrak v(x)+I\pr_\mathfrak h(x), \qquad
x\in\g,
$$
where $s>0$.

\paragraph{Algebra of $G$-invariant functions on $T^*Q$.}
The algebra $\mathcal F_2$ of $G$-invariant functions on $T^*Q$,
polynomial in momenta, can be identified with
$\mathbb{R}[\mathfrak v]^{H}$ ($\Ad_{H}$-invariant polynomials on
$\mathfrak v$). Within this identification, the Hamiltonian of the
metric $ds^2_{I}$ is given by $h(x)=\frac12\langle A x,x\rangle$,
$A=I^{-1}$, while the Hamiltonian of the normal metric $ds^2_0$ is
simply $h_0(x)=\frac12\langle x,x\rangle$. Further, the canonical
Poisson bracket on $T^* Q$ corresponds to the restriction of the
Lie-Poisson bracket (\ref{LP}) to $\mathbb{R}[\mathfrak v]^{H}$
(see Thimm \cite{Th}):
\begin{equation}
\{f,g\}_\mathfrak v(x)=-\langle x,[\nabla f(x),\nabla
g(x)]\rangle, \qquad f,g: \mathfrak v\to \R.\label{rpz}
\end{equation}

Let $\g_x$, $\mathfrak h_x$ and $\mathfrak k_{x_i}$ be isotropy
algebras of $x$ and $x_i$ in $\g$, $\mathfrak h$ and $\mathfrak
k$.

Consider the space $\mathfrak j_x \subset \mathfrak v$ spanned by
gradients of all polynomials in $\mathbb{R}[\mathfrak v]^{H}$. For
a generic point $x\in \mathfrak v$ we have (see \cite{BJ3, MP}):
\begin{eqnarray*}
&&\mathfrak j_x=([x,\mathfrak h]^\perp)\cap \mathfrak v=\{ \eta\in
\mathfrak v \, \vert\,  \langle \eta,[x,\mathfrak h]\rangle=0\}
=\{\eta\in \mathfrak v\, \vert\, [x,\eta]\subset \mathfrak v\}.\\
&&\mathfrak j_x=\{(\xi_1,\dots,\xi_n)\in\g\,\vert\, \sum_{i=1}^n
\xi_i=0, \quad \sum_{i=1}^n [x_i,\xi_i]=0\}.
\end{eqnarray*}

The Poisson bracket (\ref{rpz}) on the algebra
$\mathbb{R}[\mathfrak v]^{H}$ corresponds to the restriction of
the Lie-Poisson bivector
\begin{eqnarray}
\label{LPB} && \Lambda_x=\lambda_{x_1} \times \cdots \times
\lambda_{x_n}: \g \times \g \to \R, \\
&& \nonumber \lambda_\xi: \mathfrak k\times \mathfrak k \to \R,
\quad \lambda_\xi(\eta_1,\eta_2)=-\langle
\xi,[\eta_1,\eta_2]\rangle, \quad  \xi,\eta_1,\eta_2 \in \mathfrak
k
\end{eqnarray}
to $\mathfrak j_x$. Denote this restriction by $\bar\Lambda$. Note
that the kernel of $\bar\Lambda_x$ is \cite{BJ3}:
$$
\ker\bar\Lambda_x=\pr_\mathfrak v \ker\Lambda_x=\pr_\mathfrak
v\g_x= \pr_\mathfrak v (\mathfrak k_{x_1},\dots,\mathfrak
k_{x_n})\subset \mathfrak j_x.
$$

Here, for simplicity, the gradient operator with respect to the
restriction of $\langle\cdot,\cdot\rangle$ to $\mathfrak v$ is
also denoted by $\nabla$.

Note that \eqref{rpz} is a Poisson bracket within
$\mathbb{R}[\mathfrak v]^{H}$, while it is an almost-Poisson
bracket within the algebra of polynomials on $\mathfrak v$ (it
does not satisfies the Jacobi identity).

We have the following simple basic statement.

\begin{lemma} The differential dimension and index of
$\mathbb{R}[\mathfrak v]^{H}$ are given by:
\begin{equation}\label{DIM}
\ddim\mathbb{R}[\mathfrak v]^{H}=(n-2)\dim K, \qquad
\dind\mathbb{R}[\mathfrak v]^{H}=n\rank K.
\end{equation}
\end{lemma}

\noindent{\it Proof.} The differential dimension  of
$\mathbb{R}[\mathfrak v]^{H}$ is equal to the codimension of a
generic orbit $\Ad_H(x)$ within $\mathfrak v$, that is
$\ddim\mathbb{R}[\mathfrak v]^{H}=\dim\mathfrak j_x=\dim\mathfrak
v-\dim H+\dim\mathfrak h_x$, for a generic $x\in\mathfrak v$.
Since
\begin{equation}\label{HX}
\dim\mathfrak h_x=\dim(\mathfrak k_{x_1} \cap \dots\cap \mathfrak
k_{x_n})=0,
\end{equation}
for a generic $x\in\mathfrak v$, we obtain the first relation in
\eqref{DIM}. On the other side, from \eqref{HX} we get
\begin{eqnarray*}
\dind\mathbb{R}[\mathfrak v]^{H} &=& \dim\ker\bar\Lambda_x=\dim
\pr_\mathfrak v \g_x=\dim\g_x-\dim\mathfrak h_x\\
&=& \dim\g_x=\dim \mathfrak k_{x_1}+\dots+\dim\mathfrak k_{x_n}
=n\rank K,
\end{eqnarray*}
for a generic $x\in\mathfrak v$. $\Box$

\

Let $\R[\mathfrak k]^K$ be the algebra of $\Ad_{K}$-invariant
polynomials on $\mathfrak k$. It is generated by $r=\rank
K$-invariant homogeneous polynomials $f^1,\dots,f^r$. The algebra
of $\Ad_G$-invariant polynomials on $\mathfrak g$ is then
generated by $n\rank K$ polynomials
$$
\mathcal Z=\{f^\alpha_i=f^\alpha\circ\pr_{\mathfrak k_i}\, \vert\,
i=1,\dots,n, \, \alpha=1,\dots,r\}
$$
and the restrictions of invariants $\mathcal Z$ to $\mathfrak v$
give $n\rank K$ independent Casimir functions of
$\mathbb{R}[\mathfrak v]^{H}$.

\paragraph{Translation of Argument and Flag of Subalgebras.}
Mishchenko and Fomenko showed that the set of polynomials induced
from the invariants by shifting the argument ${\mathcal
A}=\{f^\alpha_{a,k}(x)\,\vert\, k=1,\dots,\deg f^\alpha, \,
\alpha=1,\dots,r\}$,
$$
f^\alpha_{a,t}(x)=f^\alpha(x+ta)=\sum_{k=0}^{\deg f^\alpha}
f^\alpha_{a,k}(x)t^k,
$$
is a complete commutative  set on $\mathfrak k$:
\begin{equation}\label{A}
\ddim\mathcal A=\frac12(\dim K+\rank K),
\end{equation}
for a generic $a\in \mathfrak k$ (see \cite{MF1, Bo, PY}). As was
already mentioned, the argument shift method allows us to
construct a complete commutative subalgebra in $\mathcal F_1$. Now
we shall modify the method, by using the flag of subalgebras
\eqref{flag} to use it to construct such a subalgebra within
$\mathcal F_2$.

Let
\begin{eqnarray*}
&&\mathcal B=\mathcal B_1+\mathcal B_2+\dots+\mathcal
B_{n-1}+\mathcal Z, \\
&&\mathcal B_i=\{f^\alpha_{i,k}(x)\,\vert\,k=1,\dots,\deg
f^\alpha, \,\alpha=1,\dots,r\},
\end{eqnarray*}
where polynomials $f^\alpha_{i,k}(x)$ are defined by:
\begin{equation}\label{translacija}
f^\alpha_{i,t}(x_1,\dots,x_n)=f^\alpha(x_1+\dots+x_i+tx_{i+1})=\sum_{k=0}^{\deg
f^\alpha} f^\alpha_{i,k}(x_1,\dots,x_n)t^k.
\end{equation}

\begin{theorem}\label{prva}
\item[(i)] The set $\mathcal B$ is a commutative set of
$\Ad_H$-invariant polynomials on $\g$.

\item[(ii)] The set $\mathcal B+\mu^*(\R[\mathfrak k])$ is a
complete set of polynomials on $\g$. In particular, if $\mathcal
A$ is any complete commutative set on $\mathfrak k$, then
$\mathcal B+\mu^*(\mathcal A)$ will be a complete commutative set
on $\g$.
\end{theorem}

\noindent{\it Proof.} Step 1. The polynomial in $\mathcal B$ are
$\Ad_H$-invariant. Indeed, let $h=(k,\dots,k)\in H$. Then
\begin{eqnarray*}
f^\alpha_{i,t}(\Ad_h(x)) &=&
f^\alpha(\Ad_k(x_1)+\dots+\Ad_k(x_i)+t\Ad_k(x_{i+1}))\\
&=& f^\alpha(\Ad_k(x_1+\dots+x_i+t x_{i+1}))\\
&=& f^\alpha(x_1+\dots+x_i+t x_{i+1})=f^\alpha_{i,t}(x).
\end{eqnarray*}

Step 2. The set $\mathcal B$ is involutive. Take  polynomials
$f^\alpha_{i,t_1}(x)$ and $f^\beta_{j,t_2}(x)$ given by
\eqref{translacija}.

Let $\nabla f^\alpha=\nabla f^\alpha\vert_{x_1+\dots+x_i+t_1
x_{i+1}}$ and $\nabla f^\beta=\nabla
f^\beta\vert_{x_1+\dots+x_j+t_2x_j}$. Then $[\nabla
f^\alpha,x_1+\dots+x_i+t_1x_{i+1}]=0$, $[\nabla
f^\beta,x_1+\dots+x_j+t_2x_{j+1}]=0$ and
\begin{eqnarray}
&&\nabla f^\alpha_{i,t_1}(x)=(\underbrace{\nabla
f^\alpha,\dots,\nabla f^\alpha}_i,t_1\nabla f^\alpha,0,\dots,0),
\\
&&\nabla f^\beta_{j,t_2}(x)=(\underbrace{\nabla
f^\beta,\dots,\nabla f^\beta}_j,t_2\nabla f^\beta,0,\dots,0).
\end{eqnarray}

First, consider the case $i<j$. We have
\begin{eqnarray*}
\{f^\alpha_{i,t_1},f^\beta_{j,t_2}\}(x)&=&-\langle x_1,[\nabla
f^\alpha,\nabla f^\beta]\rangle-\dots-\langle x_i,
[\nabla f^\alpha,\nabla f^\beta]\rangle\\
&&-t_1\langle x_{i+1},[\nabla f^\alpha,\nabla f^\beta]\rangle\\
&=&-\langle x_1+\dots+x_i+t_1x_{i+1},[\nabla f^\alpha,\nabla f^\beta]\rangle\\
&=& \langle [\nabla f^\alpha,x_1+\dots+x_i+t_1x_{i+1}],\nabla
f^\beta\rangle=0.
\end{eqnarray*}

Now, let $i=j$. Then we have
\begin{eqnarray*}
\{f^\alpha_{i,t_1},f^\beta_{i,t_2}\}(x)&=&-\langle
x_1+\dots+x_i,[\nabla f^\alpha,\nabla
f^\beta]\rangle \\
&&-t_1 t_2\langle x_{i+1},[\nabla f^\alpha,\nabla f^\beta]\rangle\\
&=& -\langle x_1+\dots+x_i+ t_1 x_{i+1},[\nabla f^\alpha,\nabla
f^\beta] \rangle\\
&&+t_1 \langle x_{i+1},[\nabla f^\alpha,\nabla
f^\beta]\rangle \\
&&-t_1\langle x_1+\dots+x_i+t_2x_{i+1},[\nabla f^\alpha,\nabla
f^\beta]\rangle\\
&&+
t_1 \langle x_1+\dots+x_i,[\nabla f^\alpha,\nabla f^\beta]\rangle\\
&=& t_1\langle x_1+\dots+x_{i}+x_{i+1},[\nabla f^\alpha,\nabla
f^\beta]\rangle.
\end{eqnarray*}

In the same way:
$$
\{f^\alpha_{i,t_1},f^\beta_{i,t_2}\}(x)=
 t_2\langle x_1+\dots+x_{i}+x_{i+1},[\nabla f^\alpha,\nabla
f^\beta]\rangle.
$$
Therefore $\{f^\alpha_{i,t_1},f^\beta_{i,t_2}\}=0$ for $t_1\ne
t_2$ and taking the limit $t_1 \mapsto t_2$, we get
$$
\{f^\alpha_{i,t_1},f^\beta_{i,t_2}\}=0
$$
for all $t_1, t_2$. It follows that $\{\mathcal B,\mathcal B\}=0$.
Item (i) is proved.

\

Step 3.  For a generic $x_1,\dots,x_i$, due to the
Mishchenko-Fomenko shifting of argument method, the set of
polynomials $\mathcal B_i$, considered as polynomials in variable
$x_{i+1}$, form a complete set on $\mathfrak k_{i+1}$ with respect
to the corresponding Lie-Poisson bracket. Therefore
\begin{equation}
\dim \pr_{\mathfrak k_{i+1}}\mathrm B_{i,x} x= \frac12(\dim
K+\rank K), \quad \mathrm B_{i,x}=\Span\{\nabla
f^\alpha_{i,k}(x)\}, \label{dimBi}
\end{equation}
for a generic $x\in\g$.

Let $\mathrm B_x$ be the linear space spanned by gradients of
polynomial in $\mathcal B$ at $x\in\g$. From \eqref{dimBi} we get
\begin{eqnarray}
\dim\mathrm B_x &\ge& \rank K+
\dim\pr_{\mathfrak k_2}\mathrm B_{1,x}+\dots+\dim\pr_{\mathfrak k_n}\mathrm B_{n-1,x}\nonumber\\
&\ge&  \rank K+\frac{n-1}2(\dim K+\rank K),
\label{dimB}
\end{eqnarray} where we used that $\mathcal B$
contains invariants in variable $x_1$. Thus $\ddim\mathcal
B\ge\frac12((n-1)\dim K+(n+1)\rank K)$.

\

Step 4. Fix a generic ($n\dim K-n\rank K$)-dimensional adjoint
orbit
$$
\mathcal O=\Ad_G(x_1,\dots,x_n)=\mathcal
O_1(x_1)\times\dots\times\mathcal O_n(x_n), \quad \mathcal
O_i(x_i)=\Ad_K(x_i),
$$
such that \eqref{HX} holds. This mean that the action of
$H=\diag(K)$ is locally free.

The orbit $\mathcal O$ with the Konstant-Kirillov symplectic form
$\omega$ is a symplectic leaf in $(\g,\{\cdot,\cdot\})$. The
$\Ad_H$-action, restricted to $\mathcal O$, is Hamiltonian with
the momentum mapping \eqref{mu} (e.g., see \cite{P}).

The algebra of $H$-invariant and Noether's functions
$C^\infty_H(\mathcal O)+\mu^*(C^\infty(\mathfrak k))$ is a
complete algebra on $(\mathcal O,\omega)$ and
\begin{eqnarray*}
&&\ddim C^\infty_H(\mathcal O)=(n-1)\dim K-n\rank K, \\
&&\ddim \mu^*(C^\infty(\mathfrak k))=\dim K \\
&&\ddim\left(C^\infty_H(\mathcal O)+\mu^*(C^\infty(\mathfrak k))\right)=n\dim K-(n+1)\rank K\\
 &&\dind C^\infty_H(\mathcal O)=\dind \mu^*(C^\infty(\mathfrak
k))
=\rank K\\
&&\dind\left(C^\infty_H(\mathcal O)+\mu^*(C^\infty(\mathfrak
k))\right)=\rank K
\end{eqnarray*}
 (see Theorem 2.1 and Remark 2.1 in \cite{BJ2}, where we used that a generic
$\Ad_K$-orbit in $\mu(\mathcal O)$ is regular and that the
$\Ad_H$-action is locally free at a generic point $x\in\mathcal
O$). In particular, a commutative set $\mathcal C\subset
C^\infty_H(\mathcal O)$ is a complete subset if
\begin{equation}\label{uslov}
\ddim\mathcal C=\frac12\left(\ddim C^\infty_H(\mathcal O)+\dind
C^\infty_H(\mathcal O)\right)=\frac{n-1}2(\dim K-\rank K).
\end{equation}

Let $\mathcal C=\{f\vert_\mathcal O\, \vert\, f\in \mathcal B\}$.
The invariants $\mathcal Z$ restricted to $\mathcal O$ are
constants, so we have
\begin{equation}\label{uslov*}
\ddim \mathcal C=\ddim \mathcal B-n\rank K \ge \frac12((n-1)\dim
K-(n-1)\rank K).
\end{equation}
From \eqref{uslov}, we get that $\mathcal C$ is a complete
commutative subset of $C^\infty_H(\mathcal O)$. In particular,
inequalities in \eqref{dimB} and \eqref{uslov*} are equalities.

Since the set of $\mathcal B+\mu^*(\R[\mathfrak k])$ is a complete
set restricted to a generic symplectic leaf $(\mathcal O,\omega)$,
it is a complete set on $(\g,\{\cdot,\cdot\})$. This completeness
the proof. $\Box$

\

By using Theorem \ref{prva}, we obtain the following integrable
model. Consider a left-invariant metric on $G$ defined by the
Hamiltonian function
\begin{equation}\label{novi}
\hat h_{s,t}=\sum_{i=1}^{n-1} \frac12  \langle
s_i(x_1+\dots+x_i)+t_i x_{i+1},s_i(x_1+\dots+x_i)+t_i
x_{i+1}\rangle,
\end{equation}
where parameters $s_i,t_i$ are chosen such that $\hat h$ is a
positive definite Hamiltonian of the left-invariant metric.

\begin{corollary}
The Euler equations on $\g$ determined with Hamiltonian
\eqref{novi}
\begin{eqnarray*}
&& \dot x_1=[x_1,\sum_{i=1}^{n-1}\left(s_i^2(x_1+\dots+x_i)+t_is_i
x_{i+1}\right)],\\
 &&\dot
x_k=[x_k,s_{k-1}t_{k-1}(x_1+\dots+x_{k-1})+\sum_{i=k}^{n-1}\left(s_i^2(x_1+\dots+x_i)+t_is_i
x_{i+1}\right)],\\
&&\dot x_n=[x_n,s_{n-1}t_{n-1}(x_1+\dots+x_{n-1})], \qquad \qquad
\qquad k=2,\dots,n-1
\end{eqnarray*}
are completely integrable.
\end{corollary}

\begin{lemma}{\rm \cite{BJ1}}\label{restrikcija}
\label{komutativnost} If $f$ and $g$ are $\Ad_H$-invariant
polynomials on $\g$ and $\{f,g\}=0$, then $\{f\vert_\mathfrak
v,g\vert_\mathfrak v\}_\mathfrak v=0$, where
$\{\cdot,\cdot\}_\mathfrak v$ is the bracket given by \eqref{rpz}.
\end{lemma}

Let $\mathcal F$ be the set of polynomials, obtained by
restriction of polynomials in $\mathcal B$ to $\mathfrak v$.

\begin{theorem}\label{n=3}
The set $\mathcal F$ is a complete commutative subset of
$\R[\mathfrak v]^H$.
\end{theorem}

\noindent{\it Proof.} According to Theorem \ref{prva} and Lemma
\ref{restrikcija}, the set $\mathcal F$ is commutative. Further,
from \eqref{C} and Lemma 1, it is complete if and only if
\begin{equation}
\ddim\mathcal F=\frac12((n-2)\dim K+n\rank K).
\label{dimF}\end{equation}

Since
\begin{equation*}
\mathrm F_x=\Span \{\nabla_x f(x)\, \vert \, f\in\mathcal F\}
=\Span\{\pr_\mathfrak v \nabla f(x)\,\vert\, f\in\mathcal
B\}=\pr_\mathfrak v\mathrm B_x,
\end{equation*}
we have
\begin{equation}\label{F-B}
\dim\mathrm F_x=\dim \mathrm B_x-\dim(\mathrm B_x\cap\mathfrak h).
\end{equation}

The relation \eqref{dimBi} is satisfied for a generic
$x\in\mathfrak v$ and $i<n-1$, while for  $i=n-1$ it does not
hold. Indeed, from $x_1+\dots+x_n=0$, we get that
$f^\alpha_{n-1,t}(x)=f^\alpha((1-t)x_n)=(1-t)^{\deg
f^\alpha}f^\alpha(x_n)$. Thus
\begin{eqnarray}
\dim\mathrm B_x &\ge& \rank K+ \dim\pr_{\mathfrak k_2}\mathrm
B_{1,x}+\dots+\dim\pr_{\mathfrak
k_{n-1}}\mathrm B_{n-2,x}+\rank K\nonumber\\
&=& 2\rank K+\frac{n-2}2(\dim K+\rank K),
\label{dimBv}
\end{eqnarray}
for a generic $x\in \mathfrak v$.

On the other hand, it is obvious that
\begin{equation}\label{presek}
\dim(\mathrm B_x \cap \mathfrak h) \le \dim(\Span\{\nabla
f^\alpha(x_n) \,\vert\, \alpha=1,\dots,\rank K\})=\rank K.
\end{equation}

Combining (\ref{dimBv}), \eqref{F-B} and \eqref{presek} we get
\begin{equation}
\dim \mathrm F_x \ge \frac12((n-2)\dim K +n\rank K), \label{dimF2}
\end{equation}
for a generic $x\in\mathfrak v$. According to (\ref{dimF}) we have
$\ddim\mathcal F \le\frac12((n-2)\dim K +n\rank K)$,  i.e., the
relation (\ref{dimF2}) is an equality.  $\Box$

\

Let $h_{s,t}$ be the restriction of the Hamiltonian \eqref{novi}
to $\mathfrak v$ and $ds^2_{s,t}$ be the corresponding
$G$-invariant submersion metric on $Q=G/H$.

\begin{corollary}\label{potpun}
The geodesic flow of the metric $ds^2_{s,t}$ is completely
integrable. The complete commutative set of analytic functions,
polynomial in momenta is
$$
\{\tau(f^\alpha_{i,k}\vert_\mathfrak v),
\tau(f^\alpha_i\vert_\mathfrak
v),\Phi^*(f^\alpha_{a_i,k})\,\vert\, i=1,\dots,n-1, k=1,\dots,\deg
f^\alpha, \alpha=1,\dots,r\}.
$$
\end{corollary}

Here $\tau$ denotes the bijection $\R[\mathfrak v]^H\to \mathcal
F_2$, $\Phi: T^*Q\to \g^* \cong \g$ is the momentum mapping of the
canonical $G$-action,
$$
f^\alpha(x_i+ta_i)=\sum_{k=0}^{\deg f^\alpha}
f^\alpha_{a_i,k}(x_1,\dots,x_n)t^k
$$
and $a_i\in\mathfrak k$, $i=1,\dots,n$ are in generic position.
That is, $\{f^\alpha_{a_i,k}\}$ is a complete commutative set on
$\g$ induced from the invariants by the argument translation with
$a=(a_1,\dots,a_n)$.

\paragraph{Gaudin Type Systems on $G=K^n$.}
Consider the Hamiltonian
\begin{equation*}
\hat h_a(x)=\frac12\left\langle \frac{1}{a_1} x_1 + \dots +
\frac{1}a_n x_n, \frac{1}{a_1} x_1 + \dots + \frac{1}a_n x_n
\right\rangle. \label{GH}\end{equation*}

 The
corresponding Euler equations on $\g=\mathfrak k^n$ are
\begin{equation}
\dot x_i=\sum_{j=1}^n \frac 1{a_ia_j} [x_i,x_j], \qquad
i=1,\dots,n. \label{gaudin}\end{equation}

Following \cite{P}, we refer to system (\ref{gaudin}) as a Gaudin
type system on $\mathfrak k^n$ (the Gaudin system is originally
defined for $\mathfrak k=su(2)$).

The system is $H$-invariant, so the momentum mapping (\ref{mu}) is
conserved along the flow. By using the pencil of compatible
Poisson brackets (e.g., see Bolsinov \cite{Bo}) defined by the
Lie-Poisson bivector (\ref{LPB}) and the bivector
$$
\hat{\Lambda}_x=a_1\lambda_{x_1} \times \cdots \times a_n
\lambda_{x_n}: \g \times \g \to \R,
$$
Panasyuk proved the integrability of equations (\ref{gaudin})
restricted to {\it admissible adjoint orbits} $\mathcal
O=\Ad_G(x_1,\dots,x_n)$ for a generic value of parameters
$a=(a_1,\dots,a_n)$ \cite{P}. The complete algebra of integrals is
$ \mathcal P+\mu^*(\R[\mathfrak k])$, where
\begin{equation*}
\mathcal P=\left\{f\left(\frac{x_1}{t_1 + a_1 t_2} + \dots
+\frac{x_n}{t_1 + a_n t_2}\right)\, \Big\vert\, t_1,t_2 \in \R,
t_1^2+t_2^2\ne 0, \, f\in \R[\mathfrak k]^K \right\}.
\end{equation*}

The set $\mathcal P$ is commutative. It could be proved that the
set of polynomials, obtained by the restriction of polynomials in
$\mathcal P$ to $\mathfrak v$ is a complete commutative subset of
$\R[\mathfrak v]^H$.

\section{Einstein Metrics}

Recall that the Riemannian manifold $(Q,g)$ is called {\it
Einstein} if the Ricci curvature $\mathrm{Ric}(g)$ satisfies the
equation $\mathrm{Ric}(g)=C\cdot g$, for some constant $C$
\cite{Be}.

From now on we assume  that $K$ is a simple Lie group. The normal
$G$-invariant metric $ds^2_0$ on \eqref{L-O} is Einstein (see
Proposition 5.5, \cite{WZ}). Up to the isometry and homothety, the
homogeneous space $Q$ for $n=3$ ($n\ge 4$) admits exactly
(respectively, at least) two $G$-invariant Einstein metrics that
we shall describe below.

Let $(\cdot,\cdot)_{\mathfrak v}$ be an $\ad_\mathfrak
h$-invariant scalar product on $\mathfrak v$. We can diagonalize
$(\cdot,\cdot)_{\mathfrak v}$ and $\langle
\cdot,\cdot\rangle\vert_\mathfrak v$ simultaneously (see
\cite{Nik}). Namely, let $\nu=(\nu_1,\dots,\nu_n)\in \R^n$ be an
unit vector and $\mathfrak k_\nu \subset \g$ be a linear subspace
defined by
\begin{equation}\label{modul}
\mathfrak k_\nu=\{(\nu_1 \xi,\dots,\nu_n \xi) \, \vert
\, \xi\in\mathfrak k\}.
\end{equation}
There exist $n-1$ orthogonal $\ad_\mathfrak h$-invariant
irreducible submodules $\mathfrak v_1,\cdots,\mathfrak
v_{n-1}\subset \mathfrak v$ and $n-1$ positive numbers
$s_1,\dots,s_{n-1}$ such that
$$
\mathfrak v=\mathfrak v_1 \oplus \mathfrak v_2 \oplus \dots \oplus
\mathfrak v_{n-1}, \qquad \mathfrak v_1=\mathfrak
k_{\nu^1},\dots,\mathfrak v_{n-1}=\mathfrak k_{\nu^{n-1}}
$$
and
\begin{equation}
(\cdot,\cdot)_{\mathfrak v}=s_1 \langle \cdot,\cdot
\rangle\vert_{\mathfrak v_1} \oplus s_2 \langle \cdot,\cdot
\rangle\vert_{\mathfrak v_2}\oplus \cdots \oplus s_{n-1} \langle
\cdot,\cdot \rangle\vert_{\mathfrak v_{n-1}},
\label{metric}\end{equation} where $\nu^1,\dots,\nu^{n-1}$ is the
orthonormal base of the hyperplane orthogonal to
$(1,\dots,1)\in\R^n$. The diagonalization is unique if all $s_i$
are different.

Now, let $ds^2_{p,q}$ be a $G$-invariant metric defined by the
scalar product (\ref{metric}), where
\begin{eqnarray}
&&\nu^j=\frac{1}{\sqrt{j^2+j}}(\underbrace{1,\dots,1}_j,-j,0,\dots,0),\qquad  j=1,\dots,n-1, \nonumber\\
&&s_1=\dots=s_{n-2}=1/p, \qquad s_{n-1}=1/q.
\label{s}\end{eqnarray} It is Einstein for $ p={n}^{1/(n-1)}$,
$p^{n-2}q=1$. Moreover, for $n=3$, up to isometry and homothety,
this is the only $G$-invariant Einstein metric different from the
normal one $p=q=1$ (see \cite{Nik}).

Together with the scalar product (\ref{metric}), (\ref{s}) it is
natural to consider its extension to an $\Ad_H$-invariant scalar
product on $\g$
\begin{equation}
(\cdot,\cdot)_{\mathfrak g}=\frac{1}{s}\langle \cdot,\cdot
\rangle\vert_{\mathfrak h}\oplus \frac{1}{p} \langle \cdot,\cdot
\rangle\vert_{\mathfrak v_1} \oplus \dots \oplus \frac{1}{p}
\langle \cdot,\cdot \rangle\vert_{\mathfrak v_{n-2}}\oplus
\frac{1}{q} \langle \cdot,\cdot \rangle\vert_{\mathfrak v_{n-1}}
\label{metric2}
\end{equation}
and the corresponding left-invariant metric $ds^2_{p,q,s}$ on $G$.
Then $ds^2_{p,q}$ can be seen as a submersion metric, induced by
$ds^2_{p,q,s}$.

The Hamiltonian of the metric $ds^2_{p,q,s}$, in the
left-trivialization, read
\begin{eqnarray}
\hat h &=&\frac{s}{2}  \langle \pr_\mathfrak h x,\pr_\mathfrak h x
\rangle+\frac{p}{2} \langle x-\pr_\mathfrak h x-\pr_{\mathfrak
v_{n-1}} x,x-\pr_\mathfrak h x- \pr_{\mathfrak v_{n-1}} x
\rangle\nonumber\\  &&+\frac{q}{2} \langle \pr_{\mathfrak v_{n-1}}
x,\pr_{\mathfrak v_{n-1}} x \rangle. \label{EH}
\end{eqnarray}

Note that the orthogonal projection to \eqref{modul} with respect
to the Killing form is given by
\begin{equation}\label{L3}
\pr_{\mathfrak k_\nu}(x_1,\dots,x_n)=(\nu_1(\nu_1 x_1+\dots+\nu_n
x_n),\dots,\nu_n(\nu_1 x_1+\dots+\nu_n x_n)).
\end{equation}

By using \eqref{proj-h} and \eqref{L3}, we easily get:

\begin{lemma}
The Hamiltonian \eqref{EH} has the form
\begin{eqnarray}
\hat h&=&\frac{p}{2}\sum_{k=1}^{n-1}\langle x_k,x_k\rangle +
\frac12\left(\frac{q n}{n-1}-\frac{p}{n-1}\right)\langle
x_n,x_n\rangle\nonumber \\
&&+\frac12\left(\frac{s}{n}-\frac{p}{n-1}+\frac{q}{n^2-n}\right)
\langle\mu,\mu\rangle+ \left(\frac{p}{n-1}-\frac{q}{n-1}\right)
\langle\mu,x_n\rangle \label{EH2}\end{eqnarray} where $\mu$ is the
momentum mapping \eqref{mu}.
\end{lemma}

The Euler equations with Hamiltonian (\ref{EH2}) are
\begin{eqnarray}
\nonumber && \dot x_k=[x_k, u (x_1+\dots+x_{n-1})+ vx_n  ], \qquad k=1,\dots,n-1,\\
\label{EN} && \dot x_n=[x_n, v(x_1+\dots+x_{n-1}) ],
\end{eqnarray}
where $u={s}/{n}-{p}/{(n-1)}+{q}/{(n^2-n)}$ and
$v={s}/{n}-{q}/{n}$. In particular, the set $\mathcal B$ is a set
of integral of the system (\ref{EN}). Whence, the functions
$\mathcal F$ commute with the Hamiltonian $h=\hat h\vert_\mathfrak
v$ of the metric $ds^2_{p,q}$. Applying Corollary \ref{potpun} we
obtain:

\begin{corollary}
The geodesic flow of the $G$-invariant Nikonorov's Einstein metric
on \eqref{L-O} is completely commutatively integrable by means of
analytic integrals, polynomial in momenta.
\end{corollary}

Note that, restricted to the invariant subspace $\mathfrak
v=\mu^{-1}(0)$, the equations (\ref{EN}) take the form
\begin{eqnarray}
&& \dot x_k=[x_k, (v-u)x_n  ], \qquad   k=1,\dots,n-1,\nonumber\\
&& \dot x_n=0, \label{ENR}\\
&& x_1+\dots+x_n=0.\nonumber
\end{eqnarray}
The generic solution of (\ref{ENR}) is given by
\begin{eqnarray*}
&&x_k(t)=\Ad_{\exp(t\xi)} x_k^0, \qquad x_k^0=x_k(0), \\
&&\xi=(v-u)(x_1^0+\dots +x^0_{n-1}), \quad k=1,\dots,n-1,
\end{eqnarray*}
where $\exp:\mathfrak k\to K$ is the exponential mapping.

\subsection*{Acknowledgments.} I am greatly thankful to Alexey
Bolsinov on useful discussions. This research was supported by the
Serbian Ministry of Science Project 144014 Geometry and Topology
of Manifolds and Integrable Dynamical Systems.

\end{document}